\documentclass[MSNbibl,number,citesort,dvips]{arxbj}
\usepackage{dcolumn}
\usepackage{graphicx}

\aid{0}
\volume{19}
\issue{1}
\pubyear{2013}
\firstpage{295}
\lastpage{307}
\doi{10.3150/11-BEJ393}

\makeatletter
\newtheorem{theorem}{Theorem}
\newremark{remark}{Remark}
\newcommand{\RR}{\mathbb{R}}
\newcommand{\bphi}{\bolds{\phi}}

\newcolumntype{d}[1]{D{.}{.}{#1}}
\makeatother

\begin{document}
\begin{frontmatter}

\title{Empirical likelihood-based tests for stochastic ordering}
\runtitle{Empirical likelihood and stochastic ordering}

\begin{aug}
\author[1]{\fnms{Hammou} \snm{El Barmi}\corref{}\thanksref{1}\ead[label=e1]{hammou.elbarmi@baruch.cuny.edu}}
\and
\author[2]{\fnms{Ian W.} \snm{McKeague}\thanksref{2}} %
\runauthor{H. El Barmi and I.W. McKeague} 
\address[1]{Department of Statistics and Computer Information Systems,
Baruch College, The City University of New York,
One Baruch Way, New York, NY 10010, USA}
\address[2]{Department of Biostatistics, Mailman School of Public
Health, Columbia University, 722 West 168th Street, 6th Floor,
New York, NY 10032, USA}
\end{aug}

\received{\smonth{5} \syear{2010}}
\revised{\smonth{4} \syear{2011}}

%
\begin{abstract}
This paper develops an empirical likelihood approach
to testing for the presence of stochastic ordering among univariate
distributions
based on independent random samples from each distribution. The
proposed test statistic is formed
by integrating a localized empirical likelihood statistic with respect
to the empirical
distribution of the pooled sample.
The asymptotic null distribution of this test statistic
is found to have a simple distribution-free
representation in terms of standard Brownian bridge processes.
The approach is used to compare the lengths of
rule of Roman Emperors over various historical periods, including the
``decline and fall'' phase of the empire.
In a simulation study, the power of the proposed test is found to
improve substantially upon that of a competing test due to El Barmi and
Mukerjee.
\end{abstract}

%

\end{frontmatter}

\section{Introduction}\label{sec1}

Comparing random variables in terms of their distributions
can provide an understanding of underlying causal
mechanisms and risks. In addition, knowledge of an ordering of
distributions can
be useful for increasing the efficiency of estimation procedures, as is
well documented in the literature on order restricted inference; see,
for example, the comprehensive monograph
of Silvapulle and Sen~\cite{SilSen05}.
There are many types of ordering for the comparison of univariate
distributions. These include, with increasing generality, likelihood
ratio ordering, uniform stochastic
ordering (equivalent to hazard rate ordering), stochastic ordering, and
increasing convex ordering (of interest in economics and actuarial science);
see Shaked and Shanthikumar~\cite{ShaSha06} for an overview.

The aim of this paper is to develop
an empirical likelihood approach to testing for the presence of
the classical type of stochastic ordering.
Such ordering often arises in the biomedical sciences and reliability
engineering, for example, with
lifetime distributions of human populations exposed to higher risk, or
of engineering systems under greater stress. The notion of stochastic
ordering is due to Lehmann~\cite{Leh55} who defined a random variable $X_1$
to be \textit{stochastically larger} than a random variable $X_2$ if
$F_1(x) \leq F_2(x)$ for all $x$ (with strict inequality for some $x$), where
$F_1$ and $F_2$ are the corresponding cdfs; we write this as $F_1\succ
F_2$. For
a stochastic ordering of $k$
distributions, we write $F_1\succ F_2 \succ\cdots\succ F_k$ if
$F_j(x)\le F_{j+1}(x)$ for all $x$ and $j=1,\ldots, k-1$, with strict
inequality for some $x$ and some~$j$.

There is an extensive literature on the problem of
testing for equality of two distributions against the alternative that
they are stochastically ordered. Lee and Wolf~\cite{LeeWol76} proposed a
Mann--Whitney--Wilcoxon-type test. Robertson and Wright~\cite{RobWri81} studied
the corresponding likelihood test (LRT) in the one- and two-sample cases
when the distributions are discrete. They showed that the limiting
distributions are chi-bar square. Their results indicate that, in the
two-sample case, the LRT is not asymptotically distribution free. They
also obtained the least favorable distribution in this case. Other tests
are discussed in Dykstra, Madsen and Fairbanks~\cite{DykMadFai83},
Franck~\cite{Fra84} and
Mau~\cite{Mau88}. For more than two populations, Wang~\cite{Wan96}
discussed the LRT
in the multinomial case; El Barmi and Johnson
\cite{ElBJoh06} showed that the limiting distribution of his test
statistic is of
chi-bar square type and gave the expression of the weighting values.
Also in the $k$-sample case ($k\ge2$), El Barmi and Mukerjee~\cite{ElBMuk05}
provided an asymptotically distribution-free test based on the sequential
testing procedure originally introduced by Hogg~\cite{Hog62}. This test
is applicable in
both the multinomial and the continuous cases, with or without censoring.
Recently, Baringhaus and Gr\"ubel~\cite{BarGru09} introduced a nonparametric
two-sample test for the more general hypothesis of increasing convex
ordering; their test is not
asymptotically distribution-free, however, and
requires the critical values to be obtained via a bootstrap procedure.

The contribution of the present paper
is to provide
empirical likelihood based $k$-sample tests for alternatives that are
stochastically ordered.
The empirical likelihood (EL) method was originally introduced by Owen
\cite{Owe88,Owe90} for the purpose of finding confidence regions for
parameters defined by
general classes of estimating equations. It
combines the flexibility of nonparametric methods with the efficiency
of likelihood-ratio-based inference. Inference based on EL has many
attractive properties:
estimation of variance is typically not required, the range of the
parameter space is automatically respected and confidence regions have
greater accuracy than
those based on the Wald approach. Einmahl and McKeague~\cite{EinMcK03}
developed a localized version of EL, to allow
nonparametric hypothesis testing, and showed via simulation studies that
it outperforms (in terms of power) the corresponding Cram\'er--von
Mises statistics
for a variety of classical
testing problems. Their approach is restricted to omnibus alternatives,
whereas ordered alternatives
are often more useful because they can provide a more
direct interpretation of the result of the test.

The development of the proposed test statistic
and results on its asymptotic null distribution are given in Section~\ref{sec2}.
First we consider the special case of testing whether a distribution
function is
stochastically larger than a
specified distribution function, based on a single sample. Once the
theory has been developed in this
one-sample case, it is relatively straightforward to extend the approach
to the general $k$-sample setting in which all the distribution
functions are unknown.
Section~\ref{sec3} presents the results of a simulation study in which we find
that the proposed test has
superior power to
the test of El Barmi and Mukerjee~\cite{ElBMuk05}, which is the only
previous test to have been developed
for ordered alternatives in this setting. Section~\ref{sec3} also contains
an application of the proposed test to a comparison of the lengths of
rule of Roman Emperors over various historical periods. Some concluding
remarks are given in Section~\ref{sec4}.
Proofs of the main results are collected in Section~\ref{sec5}.

\section{Empirical likelihood approach}\label{sec2}
\subsection{Stochastic ordering relative to a specified distribution}\label{sec2.1}
Suppose we are given a random sample $X_{1}, X_2, \ldots, X_n$ from the
cdf $F$,
and we want to test the null hypothesis $H_0\dvt F=F_0$ versus $ H_1\dvt F
\succ F_0$, where $F_0$ is a specified cdf.

Adapting the approach of Einmahl and McKeague~\cite{EinMcK03} to the
present setting, we first need to consider testing the ``local'' null
hypothesis $H_0^x\dvt F(x)=F_0(x)$
versus the alternative $H_1^x\dvt F(x) < F_0(x)$,
where $x$ is fixed.
The empirical likelihood procedure in this case rejects $H_0^x$ for
small values of
%
\begin{equation}
\label{el}
\mathcal{R}(x) = \frac{\sup\{ L(F)\dvt  F(x) =F_0(x)\}}{
\sup\{ L(F)\dvt  F(x) \le F_0(x)\}},
\end{equation}
where the suprema are over cdfs $F$ that are supported by the data
points, $L(F)$ is the nonparametric likelihood function and, by
convention, $\sup\emptyset=0$ and $0/0=1$.
For $F$ having point mass $p_i$ at $X_i$, define the new parameters
$\theta_{i} = p_{i}/\phi$ and
$\psi_i = p_{i}/(1-\phi)$, where
$0<\phi = F(x)<1$.
In terms of this new parameterization, with $\hat{F}$ denoting the
empirical cdf, we need to maximize
%
\begin{equation}
\label{altpar}
L(F)=\prod_{i=1}^{n}p_{i} = \biggl\{\prod_{i\dvtx  X_i\le x}\theta_{i}
\biggr\} \biggl\{\prod_{i\dvtx  X_i> x}\psi_{i}\biggr\} \phi^{n\hat{F}(x)}[1-\phi
]^{n(1-\hat{F}(x))},
\end{equation}
subject to the constraint
\[
\sum_{i\dvtx  X_i\le x} \theta_{i} = \sum_{i\dvtx  X_i > x} \psi_{i} = 1,
\]
with either $\phi= F_0(x)$ under ${H}_0^x$, or $\phi< F_0(x)$ under ${H}_1^x$.
Note that the three terms in the right-hand side of (\ref{altpar}) can
be maximized separately.
As the constraints for the first two terms of (\ref{altpar}) are the
same for both the numerator and the denominator of (\ref{el}),
these terms cancel and make no contribution to $\mathcal{R}(x)$.
The third term of (\ref{altpar}) is maximized by
$\phi= F_0(x)$ under ${H}_0^x$, or $ \phi= F_0(x) \wedge\hat{F}(x)$
under ${H}_1^x$.
Consequently,
\[
\mathcal{R}(x)  =  \cases{
1, & \quad$\mbox{if }  \hat{F}(x) > F_0(x),$\vspace*{2pt}\cr
\biggl[\displaystyle\frac{F_0(x)}{\hat F(x)}\biggr]^{n\hat{F}(x)}\biggl[\displaystyle\frac
{1-F_0(x)}{ 1-{\hat F}(x)} \biggr]^{n(1-\hat{F}(x))} ,&\quad$\mbox{if }\hat{F}(x)\le
F_0(x),$}
\]
with the convention that any term raised to a zero power is set to 1.
Using a second-order Taylor expansion of $\log(1+y)$ about $y=0$, it
can be shown (see the proof of the theorem below) that, for
a given $x$, such that $0<F_0(x)<1$, under $H_0^x$,
\begin{eqnarray*}
\label{Taylor}
-2\log\mathcal{R}(x) &=& n\bigl(\hat{F}(x)- F_0(x)\bigr)^2\biggl[\frac{1}{\hat{F}(x)}
+ \frac{1}{1-\hat{F}(x)}\biggr] I[0<\hat{F}(x)\le F_0(x)] +\mathrm{o}_p(1) \\
& \stackrel{d}{\rightarrow} & Z^2 I(Z\ge0),
\end{eqnarray*}
using the CLT and the continuous mapping theorem, where $Z\sim N(0,1)$.
That is, the asymptotic null distribution of $-2\log\mathcal{R}(x)$ is
chi-bar square.

To test ${H}_0$ against ${H}_1$, we introduce the integral-type test statistic
\[
T_n = -2\int_{-\infty}^{\infty}\log(\mathcal{R}(x))  \,\mathrm{d}F_0(x).
\]
Here the range
of integration is actually restricted to the interval
$[{X_{(1)}},{X_{(n)}}]$, where $X_{(1)}$ and $X_{(n)}$ are the smallest
and largest order statistics in the sample, because the integrand
vanishes outside this interval. The following result
gives the asymptotic null distribution of $T_n$.
\begin{theorem}
\label{th1}
If $F_0$ is continuous, then under $H_0,$
\[
T_n \stackrel{d}{\rightarrow} \int_0^1 \frac{B^2(t)}{t(1-t)}I\bigl(B(t) \ge
0\bigr)  \,\mathrm{d}t,
\]
where $B$ is a standard Brownian bridge.
\end{theorem}

\begin{remark}\label{rem1} An alternative test statistic is obtained
by integrating with respect to the empirical cdf (instead of $F_0$),
\[
T_n^*= -2\int_{-\infty}^{\infty}\log(\mathcal{R}(x))  \,\mathrm{d}\hat F(x).
\]
It can be shown using a martingale argument (see Section~\ref{sec5}), that
$T_n^*$ has the same asymptotic null distribution as $T_n$.
\end{remark}

\subsection{Stochastic ordering among $k$ distributions}\label{sec2.2}

Suppose now that we are given a random sample of size $n_j$
from the cdf $F_j$, for $j=1,\ldots, k$,
the $k$ samples are independent and we want to test the null hypothesis
$H_0\dvt  F_1=\cdots=F_k$ versus
$H_1\dvt  F_1\succ \cdots\succ F_k$.
We assume that
the proportion $w_j=n_j/n$ of observations in the $j$th sample remains
fixed as the total sample size $n\to\infty$, with $0<w_j<1$ for all
$j=1,\ldots,k$.

Adapting the approach of Section~\ref{sec2.1}, we now consider the localized
empirical likelihood function
%
\begin{equation}
\label{ksamp}
\mathcal{R}(x) = \frac{\sup\{\prod_{j=1}^k L(F_j)\dvt
F_j(x)=F_{j+1}(x), j=1,\ldots,k-1\}}
{\sup\{\prod_{j=1}^k L(F_j)\dvt  F_j(x)\le F_{j+1}(x), j=1,\ldots
,k-1\}},
\end{equation}
where, in each supremum, $F_j$ is supported by the
observations in the $j$th sample.
Applying the same parameterization used in (\ref{altpar}), separately
for each $F_j$, and making
the same cancelation in the numerator and denominator, it suffices to maximize
%
\begin{equation}
\label{altpar2}
\prod_{j=1}^k\phi_j^{n_j\hat{F}_j(x)}[1-\phi_j]^{n_j(1-\hat{F}_j(x))}
\end{equation}
subject to the constraint
$0<\phi_1=\cdots=\phi_k<1$, or $0<\phi_1\le\cdots\le\phi_k<1$,
depending on whether it is the numerator or the denominator of (\ref{ksamp}).
Here $\hat F_j$ is the empirical cdf based on the $j$th sample.
Under the first of these constraints, (\ref{altpar2})
is maximized by $\phi_j=\hat F(x)$, where $\hat F$ is the empirical
cdf of the pooled sample.
Under the second constraint, this is the classical
bioassay problem, as discussed in Robertson \textit{et al.}~\cite{RobWriDyk88}, page 32,
and it follows that (\ref{altpar2})
is maximized by
\[
\phi_j=
E_\mathbf{w}(\hat{\bphi}|\mathcal{I})_j\equiv\tilde F_j(x),
\]
where $ E_\mathbf{w}({\hat{\bphi}}|\mathcal{I})$ is the weighted least squares
projection of ${\hat{\bphi}}=(\hat F_1(x),\ldots,\hat F_k(x))^T$ onto
$\mathcal{I} = \{\mathbf{z} \in
\RR^k\dvt    z_1 \leq z_2 \leq\cdots\leq z_k\}$, with weights $w_j$. In
passing, we mention that several algorithms have been developed for
computing this projection, including the pool-adjacent-violators
algorithm, see Robertson \textit{et al.}~\cite{RobWriDyk88}.
We now have
%
\begin{equation}
\label{likratio}
\mathcal{R}(x) = \prod_{j=1}^k \biggl[\frac{\hat F(x)}{\tilde F_j(x)}
\biggr]^{n_j\hat{F}_j(x)}\biggl[\frac{1-\hat F(x)}{ 1-{\tilde F}_j(x)}
\biggr]^{n_j(1-\hat{F}_j(x))}
\end{equation}
under the convention that any term raised to a zero power is set to 1.

To test $H_0$ against $H_1$, we propose the test statistic
%
\begin{equation}
{T}_n = -2\int_{-\infty}^{\infty} \log\mathcal{R}(x)  \,\mathrm{d}\hat F(x).
\end{equation}
The following theorem gives the asymptotic null distribution of $T_n$.
\begin{theorem}
\label{th2}
Under $H_0$ and assuming that the common distribution function $F$ is
continuous,
%
\begin{equation}
T_n \stackrel{d}{\rightarrow} \sum_{j=1}^k w_{j}\int_0^1
\frac{(E_\mathbf{w}[\mathbf{B}(t)|\mathcal{I}]_j-\overline{B}(t))^2}{t(1-t)}  \,\mathrm{d}t,
\label{limit}
\end{equation}
where $\mathbf{B}= (B_1/\sqrt{w_1},B_2/\sqrt{w_2}, \ldots, B_k/\sqrt
{w_k})^T$, the processes $B_1, B_2, \ldots, B_k$
are independent standard Brownian bridges, and $\overline{B}= \sum
_{j=1}^k \sqrt w_j B_j.$
\end{theorem}

\begin{remark}\label{rem2} For the two-sample case, it can be shown
that the limiting distribution in the above result
coincides with that in the one-sample case (Theorem~\ref{th1});
the equivalence arises from the fact that $B=\sqrt{w_2} B_1 - \sqrt
{w_1} B_2$
is a standard Brownian bridge. Moreover, when testing against the
unrestricted alternative $F_1\neq F_2$,
the limiting
distribution of the corresponding test statistic (see Einmahl and
McKeague~\cite{EinMcK03}, Theorem 2a) is the same apart from the
presence of the indicator $I(B(t) \ge0)$
in the integrand.
\end{remark}

%
\begin{table}
\tablewidth=200pt
\caption{Selected critical points of $T_n$}\label{tab1}
\begin{tabular*}{200pt}{@{\extracolsep{\fill}}llll@{}}
\hline
& \multicolumn{3}{l@{}}{Significance level $\alpha$}\\[-6pt]
& \multicolumn{3}{l@{}}{\hrulefill}\\
$k$ & 0.01 & 0.05 & 0.10 \\
\hline
2 & 3.185 & 1.821 & 1.288\\
3 & 4.128 & 2.613 & 1.943 \\
4 & 4.663 & 3.107 & 2.404\\
5 & 5.144 & 3.470 & 2.701\\
\hline
\end{tabular*}
\end{table}

\section{Numerical examples}\label{sec3}

In this section we discuss some numerical examples illustrating the
proposed test
for a comparison of two or more distributions developed in Section~\ref{sec2.2}.

To implement the proposed test we first need to obtain critical values
for $T_n$. The null distribution of $T_n$ is not tractable, even
asymptotically, but it
is asymptotically distribution free. We use simulation to approximate
selected critical values
as provided in Table~\ref{tab1}. These critical values
are based on 100\,000 data sets distributed as $N(0,1)$, with sample
sizes of $n_i=100$, $i=1,\ldots,k$, in each case.
The (Fortran) program used to compute the critical values in Table~\ref{tab1}
is available online in the supplemental files.\looseness=1

\subsection{Simulation study}\label{sec3.1}

Here we present the results of a simulation study designed to compare
the performance of
$T_n$ with the test statistic $S_n$ of El Barmi and Mukerjee \cite
{ElBMuk05}, which is defined
as the maximum of a sequence of (one-sided) two-sample
Kolmogorov--Smirnov test statistics.
As far as we know, $S_n$ is the only previously developed test
statistic when $k\ge3$.

Tables~\ref{tab2} and~\ref{tab3} give the results for a variety of
distributions and sample sizes,
for $k=2$ and $k=3$, respectively. In each case, 10\,000 data sets were
used to
approximate the power at a nominal level of $\alpha=0.05$, with
critical values for $T_n$ taken from Table~\ref{tab1};
critical values for $S_n$ are obtained from its asymptotic
distribution, which
is available in a closed form. In all cases, $T_n$ has greater power
than $S_n$ and has
better agreement with the nominal level of the test.

%
\begin{table}
\caption{Power comparison of tests for stochastic ordering of $k=2$
distributions at level \mbox{$\alpha= 0.05$}}
\label{tab2}
\begin{tabular*}{\textwidth}{@{\extracolsep{4in minus 4in}}llllllll@{}}
\hline
\multicolumn{2}{@{}l}{Distributions} & \multicolumn{2}{l}{$n_1=50$, $n_2=30$}& \multicolumn{2}{l}{$n_1= 30$, $n_2=50$}
&\multicolumn{2}{l@{}}{$n_1=50$, $n_2=50$}\\[-6pt]
\multicolumn{2}{@{}l}{\hrulefill} & \multicolumn{2}{l}{\hrulefill}& \multicolumn{2}{l}{\hrulefill}
&\multicolumn{2}{l@{}}{\hrulefill}\\
$F_1$ & $F_2$ & $T_n$ & $S_n$ & $T_n$ & $S_n$ & $T_n$ & $S_n$\\
\hline
Uni($0,1$) & Uni($0,1$) & 0.064 & 0.038 & 0.051 & 0.045 & 0.051 & 0.036\\
Uni($0,1.1$)& Uni($0,1$) & 0.143 & 0.104 & 0.162 & 0.111 & 0.199 & 0.125\\
Uni($0,2$) & Uni($0,1$) & 0.911 & 0.816 & 0.912 & 0.818 & 0.908 & 0.815\\
Uni($0.1,1.1$)&Uni($0,1$) & 0.377 & 0.244 & 0.357& 0.246 & 0.468 &
0.287\\[5pt]
Exp(1) & Exp(1) & 0.063 & 0.037 & 0.048 & 0.041 & 0.047 & 0.036\\
Exp(1) & Exp(1.1) & 0.123 & 0.076 & 0.091 & 0.068 & 0.108 & 0.076\\
Exp(1) & Exp(2) & 0.782 & 0.716 & 0.813 & 0.718 & 0.909 & 0.815\\
0.1${}+{}$Exp(1) & Exp(1) & 0.207 & 0.118 & 0.137 & 0.105 & 0.195 & 0.127\\[3pt]
$N(0,1)$ & $N(0,1)$ & 0.063 & 0.037 & 0.049 & 0.040 & 0.051 & 0.036\\
$N(0.1,1)$ & $N(0,1)$ & 0.132 & 0.081 & 0.100 & 0.079 & 0.122 & 0.079\\
$N(0.5,1)$ & $N(0,1)$ & 0.646 & 0.530 & 0.690 & 0.540 & 0.771 & 0.628 \\
$N(1,1)$ & $N(0,1)$ & 0.992 & 0.975 & 0.991 & 0.975 & 0.993 & 0.976 \\
\hline
\end{tabular*}
\end{table}
%

%
\begin{table}[b]
\caption{Power comparison of tests for stochastic ordering of $k=3$
distributions at level \mbox{$\alpha= 0.05$}}
\label{tab3}
\begin{tabular*}{\textwidth}{@{\extracolsep{4in minus 4in}}lllllll@{}}
\hline
\multicolumn{3}{@{}l}{Distributions} & \multicolumn{2}{l}{$n_1=n_2=n_3=30$}& \multicolumn{2}{l@{}}{$n_1=n_2=n_3=50$}
\\[-6pt]
\multicolumn{3}{@{}l}{\hrulefill} & \multicolumn{2}{l}{\hrulefill}& \multicolumn{2}{l@{}}{\hrulefill} \\
$F_1$ & $F_2$ & $F_3$ & $T_n$ & $S_n$ & $T_n$ & $S_n$ \\
\hline
Uni($0,1$) & Uni($0,1$) & Uni($0,1$) & 0.038 & 0.033 & 0.045 & 0.039\\
Uni($0,1.1$) & Uni($0,1$)& Uni($0,1$) & 0.455 & 0.370 & 0.740 & 0.647\\
Uni($0,1.1$) & Uni($0,1.1$)& Uni($0,1$) & 0.389 & 0.319 & 0.651 & 0.633\\
Uni($0.1,1.1$)& Uni($0,1$) & Uni($0,1$) & 0.948 &0.884 & 0.999 &0.885\\[5pt]
Exp(1) & Exp(1) & Exp(1) & 0.041 &0.019 & 0.049 & 0.045 \\
Exp(1) & Exp(1) & Exp(1.1) & 0.076 &0.033 & 0.098 & 0.067\\
Exp(1) &Exp(1.1) & Exp(1.1) & 0.067 & 0.029 & 0.098 & 0.073\\
Exp(1) & Exp(1.1) & Exp(1.2) & 0.116 & 0.046 & 0.171 & 0.109\\
Exp(1) & Exp(1.25)& Exp(1.5) & 0.313 & 0.121 & 0.507 & 0.321\\[5pt]
$N(0, 1)$ & $N(0,1)$ & $N(0,1)$ & 0.042 & 0.035 & 0.049 & 0.035\\
$N(0.1,1)$ & $N(0,1)$ & $N(0,1)$ & 0.272 & 0.183 & 0.423 & 0.292\\
$N(0.1,1)$ & $N(0.1,1)$ & $N(0,1)$ & 0.246 & 0.151 & 0.393 & 0.249\\
$N(0.5, 1)$ & $N(0.25,1)$ & $N(0,1)$ &1.000 & 0.993 & 1.000 & 1.000 \\
\hline
\end{tabular*}
\end{table}
%

\subsection{Lengths of rule of Roman Emperors}\label{sec3.2}

A recent article of Khmaladze, Brownrigg and Haywood~\cite{KhmBroHay07}
reached the interesting conclusion that the lengths of rule of Roman
Emperors were exponentially distributed, implying that their reigns
ceased unexpectedly
(``brittle power''). It is also of interest to examine whether there
were \textit{changes} in
the distribution of rule lengths, especially during the ``decline and
fall'' phase of the empire. We use the list of $n=70$ Roman Emperors
from Augustus to Theodossius, covering 27 BC to 395 AD. Our analysis
is based on the chronology of Parkin (see Khmaladze \textit{et al.} \cite
{KhmBroHay07} for further details).
The (Fortran) programs used for the two analyzes are available online
in the supplemental files.

\begin{figure}

\includegraphics{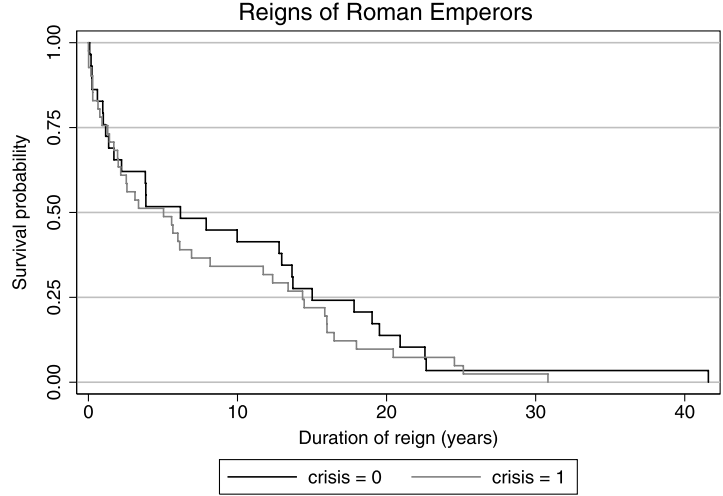}

\caption{Empirical survival functions of durations of rule of the
first 70 Roman Emperors before 235 AD (crisis${} = {}$0), and after 235 AD
(crisis${} = {}$1).}
\label{roman-emps1}
\end{figure}

First we consider whether there is an effect on duration of rule due
to the
Crisis of the Third Century (235--284 AD), when the Roman Empire nearly
collapsed under the pressure of civil war (among other things!). Figure
\ref{roman-emps1} shows the empirical survival function of durations of
rule for
the Principate (27 BC--235 AD), which was the relatively stable period
preceding the Crisis, compared with
the period after 235 AD; the sample sizes are $n_1=29$ and $n_2=41$,
respectively.
The two distributions appear to be exponential, and
the likelihood ratio test of stochastic ordering under this assumption
has $p$-value $0.195$; the corresponding unrestricted likelihood-ratio
test has $p$-value
$0.390$. Applying our proposed test (with $k=2$) to assess whether the
duration of rule is stochastically shorter
after the Principate, we obtain $T_n=0.3161$ with a $p$-value of $0.424$.
This compares with
a $p$-value of 0.575 based on $S_n$.

%

\begin{figure}

\includegraphics{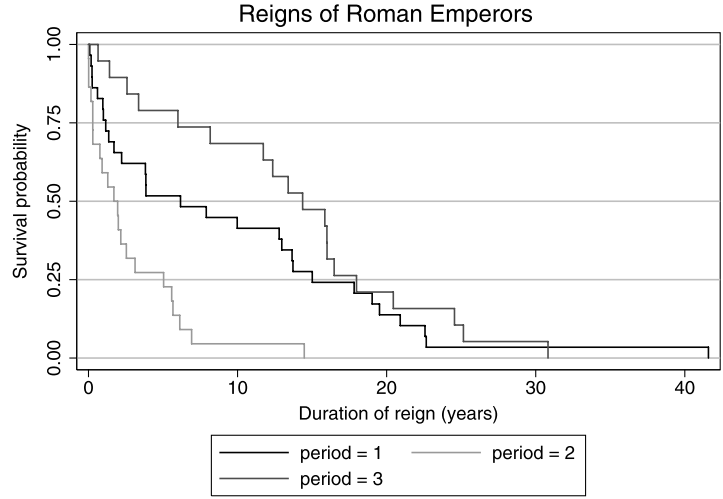}

\caption{Empirical survival functions of durations of rule during the
Principate, 27 BC--235 AD (period${} = {}$1), the
Crisis, 235--284 AD (period${} = {}$2) and the Dominate, 284--395 AD (period${} = {}$3).}
\label{roman-emps2}\vspace*{-3pt}
\end{figure}

The period 285--395 AD forms part of what is known as the Dominate, the
despotic later phase of the empire. Inspection of Figure \ref
{roman-emps2} suggests that
the exponential hypothesis is not tenable for each separate period, so
our nonparametric approach is more reasonable.
The plot also suggests that
the rule lengths are stochastically ordered as Dominate $\succ$
Principate $\succ$ Crisis. Applying our approach to formally test this
hypothesis, we find that
$T_n$ has a $p$-value of $ 0.0002$, compared with a $p$-value of 0.0017 for
$S_n$. Under the assumption of
exponential distributions, the likelihood ratio test has $p$-value less
than $0.0007$.\vspace*{-3pt}

\section{Discussion}\label{sec4}\vspace*{-3pt}
In this paper we have developed a novel empirical likelihood approach to
the important problem of nonparametrically testing for the presence of
stochastic ordering
based on $k$ independent samples.
The proposed tests are computationally efficient to implement,
and could be used with massive data sets because they do not rely on
the bootstrap
or any other simulation technique, and they reduce to a local test for
an ordering of binomial probabilities, which
only requires a single sweep through the pooled data in the $k$ groups.

%

Various extensions of the proposed tests are possible.
In change-point problems, for example, it is of interest to test
whether there is a sudden change
in the distribution of a sequence of independent random variables
$X_1,\ldots,X_n$.
Einmahl and McKeague~\cite{EinMcK03} developed an EL-based change-point test
for the presence of an (unknown) change-point $\tau\in\{2,\ldots,n\}$
such that
\[
X_1,\ldots,X_{\tau-1}\sim F_1 \quad \mbox{and} \quad   X_\tau,\ldots,X_n\sim F_2.
\]
They only considered the unrestricted alternative $F_1\neq F_2$, but
it is also of interest to consider the ordered alternative\vadjust{\goodbreak} $F_1\succ F_2$.
This can be done by extending the two-sample case to allow the sample
sizes to depend
on an additional local parameter, namely $t\in[1/n,1)$, with
$n_1=\lfloor nt\rfloor$ and $n_2=n-\lfloor nt\rfloor$.
The resulting test statistic has a limiting distribution of the same form
as in Theorem 2 of Einmahl and McKeague~\cite{EinMcK03}, involving the
integral of a four-sided tied-down
Wiener process $W_0(t,y) $, except that the integrand now includes the
indicator $I(W_0(t,y) \ge0)$.

Our approach also naturally extends to non-monotonic alternatives,
namely to
testing whether $F_1, F_2 , \ldots, F_k $ are isotonic with respect
to a quasi-order on $\{1, 2, \ldots, k\}$. A~relation $\lesssim$ on $\{
1, 2, \ldots, k\}$ is a \textit{quasi-order} if it is
reflexive and transitive (and a partial order if, in addition, it is
antisymmetric). We say that $F_1, F_2 , \ldots, F_k $ are \textit{isotonic}
with respect to $\lesssim$ if $F_i \succ F_j$ whenever $i\lesssim j$.
Examples of such ordered alternatives include $ F_1 \succ F_i$,
$i=2,\ldots, k$ (tree
ordering) and $ F_1 \succ F_2 \succ\cdots\succ F_{i_0} \prec F_{i_0+1}\prec
\cdots\prec F_k$, where $i_0$ is known (umbrella ordering).
The localized empirical likelihood (\ref{ksamp}) extends naturally to
such ordered alternatives,
the only difference being that in $\phi_j=
E_\mathbf{w}(\hat{\bphi}|\mathcal{I})_j$ the set $\mathcal{I}$ is now the
isotonic cone corresponding to $\lesssim$.
For example, in the case of tree
ordering, the cone becomes $\mathcal{I} = \{\mathbf{z} \in
\RR^k\dvt  z_1 \leq z_i, i=2,  \ldots, k\}$.
The $\phi_j$ can be computed using quadratic programming or algorithms
described in Robertson, Wright
and Dykstra~\cite{RobWriDyk88}, one of the most general being the
lower-sets algorithm. The limiting distribution of the resulting test
statistic is
obtained by taking $\mathcal{I}$ in (\ref{limit}) as the isotonic cone
corresponding to $\lesssim$.

An important and challenging problem for future research in this area would
be to develop EL-based tests for stochastic ordering based on {\it
censored} data.
EL methods are well developed for the comparison of survival functions
from right-censored data, see McKeague and Zhao \cite
{McKZha02,McKZha05}, but these methods only apply to
omnibus alternatives. The complication in extending the present tests
to right-censored data arises because the EL ratio would then no longer
have such an
explicit form as in (\ref{likratio}), and Lagrange multipliers
would be involved. This extension is beyond the scope of the present
paper.\vspace*{-3pt}

\section{Proofs}\label{sec5}\vspace*{-3pt}

\begin{pf*}{Proof of Theorem \protect\ref{th1}}
For $0<\varepsilon<1$, let
$x_{\varepsilon},y_\varepsilon$ be real numbers such that $F_0(x_{\varepsilon})
= 1-F_0(y_{\varepsilon})=\varepsilon/2$. Then decompose the test statistic as
$T_n = T_{1n}+ T_{2n}$, where
\[
T_{1n}= -2 \int_{x_{\varepsilon}}^{y_{\varepsilon}} \log(\mathcal{R}(x))   \,\mathrm{d}F_0(x)\vspace*{-3pt}
\]
and
\[
T_{2n}= -2 \int_{[x_{\varepsilon},y_\varepsilon]^c} \log(\mathcal{R}(x))  \,\mathrm{d}F_0(x).
\]
By appealing to Theorem 4.2 of Billingsley~\cite{Bil68}, note that, to
complete the proof of the theorem, it suffices to show that
for fixed $\varepsilon$,
%
\begin{equation}
\label{T1n}
T_{1n} \stackrel{d}{\rightarrow} \int_{\varepsilon/2}^{1-\varepsilon/2} \frac
{B^2(t)}{t(1-t)}I\bigl(B(t)\ge0\bigr)  \,\mathrm{d}t
\end{equation}
as $n\to\infty$, and, for each $\delta>0$, that
$ \limsup_{n\to\infty} P(|T_{2n} |\ge\delta) \to0$
as $\varepsilon\to0$.\vadjust{\goodbreak}

First consider $T_{1n}$. Using the inequality $|\log(1+y)-y+y^2/2|\le
|y|^3/3$ when $|y|\le1/2$, the Glivenko--Cantelli theorem and
Donsker's theorem,
we have
\begin{eqnarray*}
&&\limsup_{n\to\infty}\sup_{x\in[x_\varepsilon, y_\varepsilon]}\biggl|\log
(\mathcal{R}(x))+{n\over2} \bigl(\hat{F}(x)- F_0(x)\bigr)^2\biggl[\frac{1}{\hat
{F}(x)} + \frac{1}{1-\hat{F}(x)}\biggr] I[\hat{F}(x)\le F_0(x)]\biggr|\\
&&\quad   \le\limsup_{n\to\infty}\sup_{x\in[x_\varepsilon,
y_\varepsilon]} {n\over3} |\hat{F}(x)- F_0(x)|^3\biggl[\frac{1}{\hat
{F}(x)} + \frac{1}{1-\hat{F}(x)}\biggr] = 0,
\end{eqnarray*}
almost surely. Then, noting that $\hat F(x)=\hat\Gamma(F_0(x))$, where
$\hat\Gamma$ is the empirical cdf of $V_i=F_0(X_i)\sim U(0,1)$,
$i=1,\ldots, n$, and changing variables in the integration to
$t=F_0(x)$, it follows that
%
\begin{eqnarray}\label{lead}
T_{1n}&=&\int_{\varepsilon/2}^{1-\varepsilon/2} n\bigl(\hat{\Gamma}(t)- t\bigr)^2
\biggl[\frac{1}{\hat{\Gamma}(t)} + \frac{1}{1-\hat{\Gamma}(t)}\biggr] I\bigl[
\sqrt n\bigl(\hat\Gamma(t)-t\bigr)\le0\bigr]   \,\mathrm{d}t +\mathrm{o}_p(1)
\nonumber
\\[-8pt]
\\[-8pt]
\nonumber
&=& \int_{\varepsilon/2}^{1-\varepsilon/2} \frac{\hat U(t)^2}{t(1-t)} I[\hat
U(t)\le0]  \,\mathrm{d}t +\mathrm{o}_p(1),
\end{eqnarray}
where $\hat U(t)=\sqrt n(\hat\Gamma(t)-t)$ is the uniform empirical process.
Note that (for any fixed $0<\varepsilon<1$) the functional
\[
f\mapsto\int_{\varepsilon/2}^{1-\varepsilon/2} \frac{f(t)^2}{t(1-t)} I\bigl(f(t)
\le0\bigr)   \,\mathrm{d}t, \qquad f\in D[0,1],
\]
is continuous when the Skorohod space $D[0,1]$ is equipped with the
uniform norm.
By Donsker's theorem, $\hat U$ converges weakly to $B$ in $D[0,1]$, so
applying the continuous mapping theorem to the leading term in (\ref{lead}) establishes
(\ref{T1n}).

Finally we need to verify the claim concerning $T_{2n}$. This follows
immediately from a
corresponding result in Einmahl and McKeague~\cite{EinMcK03}, who
considered the test
of the null hypothesis $F=F_0$ versus the (omnibus) alternative $F\neq
F_0$, with the same integral-type test statistic
as $T_n$ except that the integrand does not vanish when $\hat F(x)>F_0(x)$.
This completes the proof.
\end{pf*}

\begin{pf*}{Proof for Remark \protect\ref{rem1}}
The asymptotic distribution of
$T_n^*$ can be obtained following the same steps
as the proof of Theorem~\ref{th1}, except that the leading term in
$T_{1n}$ now becomes
\[
\int_{\varepsilon/2}^{1-\varepsilon/2}  \frac{\hat U(t)^2}{t(1-t)} I[\hat
U(t)\le0]   \,\mathrm{d}\hat\Gamma(t)=\int_{\varepsilon/2}^{1-\varepsilon/2}   V(t-)
\,\mathrm{d}\hat\Gamma(t) +\mathrm{o}_p(1),
\]
where
\[
V(t)=\frac{\hat U(t)^2}{t(1-t)} I[\hat U(t)\le0, \varepsilon/2 < t \le
1-\varepsilon/2].
\]
Note that
\[
M(t)=\hat\Gamma(t)- \int_0^t [1-\hat\Gamma(s-)](1-s)^{-1}  \,\mathrm{d}s
\]
is a martingale wrt to the natural filtration defined by $\hat\Gamma$,
and its
predictable quadratic variation process is
$\langle M \rangle(t) =n^{-1} \int_0^t [1-\hat\Gamma(s-)](1-s)^{-1}
\,\mathrm{d}s $. Also note that $V(t-)$
is a predictable process because it is adapted and left-continuous. Write
\begin{eqnarray*}
&&\int_{\varepsilon/2}^{1-\varepsilon/2}   V(t-)   \,\mathrm{d}\hat\Gamma(t)\\
&& \quad =\int_{\varepsilon/2}^{1-\varepsilon/2}   V(t-)   \,\mathrm{d}M(t)+
\int_{\varepsilon/2}^{1-\varepsilon/2}   V(t-)[1-\hat\Gamma(t-)](1-t)^{-1}  \,\mathrm{d}t.
\end{eqnarray*}
Using a basic property of martingale integrals, the second moment of
the first term above
is
\[
E\int_{\varepsilon/2}^{1-\varepsilon/2}   V(t-)^2 \, \mathrm{d}\langle M \rangle(t)=\mathrm{O}(1/n),
\]
so this term tends in probability to zero.
The second term in the above display can be handled in the same way
as the main term $T_{1n}$ in the proof of Theorem~\ref{th1}, and has
the same limit distribution.
\end{pf*}

\begin{pf*}{Proof of Theorem \protect\ref{th2}} The proof is similar to the
proof of Theorem~\ref{th1}, so we only
indicate the main steps. Using the Taylor expansion of $\log(1+y)$, as
before, and the (uniform) consistency of $\tilde F_j$ as an estimator
of $F_j=F$ (see, e.g., El Barmi and Mukerjee~\cite{ElBMuk05}, page~253),
for each fixed $x$, such that $0<t=F(x)<1$, we have
\begin{eqnarray*}
\label{Taylor-k}
-2\log\mathcal{R}(x) &=& \sum_{j=1}^k n_j\bigl(\tilde F_j(x)- \hat{F}(x)\bigr)^2
\biggl[\frac{1}{\tilde{F}_j(x)} + \frac{1}{1-\tilde{F}_j(x)}\biggr] +\mathrm{o}_p(1) \\
& =& \sum_{j=1}^k w_j{[\sqrt n(\tilde F_j(x)- F(x)) -\sqrt n (\hat
{F}(x)-F(x))]^2\over F(x)(1-F(x))} +\mathrm{o}_p(1)\\
& =& \sum_{j=1}^k w_j\frac{(E_\mathbf{w}[\hat\mathbf{U}(t)|\mathcal{I}]_j-\overline{U}(t))^2}{t(1-t)} +\mathrm{o}_p(1)\\
&\stackrel{d}{\rightarrow} &\sum_{j=1}^k w_{j}
\frac{(E_\mathbf{w}[\mathbf{B}(t)|\mathcal{I}]_j-\overline{B}(t))^2}{t(1-t)},
\end{eqnarray*}
where $\hat\mathbf{U}=(\hat U_1/\sqrt{w_1},\hat U_2/\sqrt{w_2}, \ldots,
\hat U_k/\sqrt{w_k})^T$, $\hat U_j(t)=\sqrt{n_j}(\hat F_j(x)-F(x))$
are independent uniform empirical processes, and
$\overline{U}= \sum_{j=1}^k \sqrt w_j \hat U_j.$ Donsker's theorem and
the continuous mapping theorem
have been used as before, but we have also used the fact that $ E_{\bf
w}(\cdot|\mathcal{I})$ is a continuous
function on $\RR^k$.
\end{pf*}

\section*{Acknowledgements} The authors thank Estate Khmaladze for
sending the data on the Roman Emperors and a referee and an associate
editor for their helpful comments that have a resulted in a much
improved paper. The work of Hammou El Barmi was supported by PSC-CUNY
Grant 62795-00 40 and the work of Ian McKeague was supported in part by
NSF Grant DMS-08-06088 and NIH Grant R01 GM095722.

\begin{supplement}[id=suppA]
\sname{Supplement}
\slink[doi]{10.3150/11-BEJ393SUPP} 
\sdatatype{.zip}
\sfilename{bej393\_supp.zip}
\sdescription{We provide the (Fortran) programs as well
as the data used in the Roman Emperors example,
and the program used to compute the critical values in Table~\ref{tab1}.}
\end{supplement}

%

\printhistory

\end{document}